\documentclass[12pt]{article}%
\usepackage{amssymb}
\usepackage{amsmath}
\usepackage{amsfonts}
\usepackage{graphicx}%
\begin{document}

\title{Twin Primes and a Primality Test by Indivisibility }
\author{M. Chaves\\Escuela de Fisica, Universidad de Costa Rica\\e-mail: mchaves@fisica.ucr.ac.cr}
\date{June 20, 2009}
\maketitle

\begin{abstract}
In Wilson's Theorem the primality of a number hinges on a congruence. We
present a similar test where the primality of a number $m$ hinges, instead, on
the indivisibility of $4(m-5)!$ by $m$. One implication of this theorem is a
necessary and sufficient condition for two numbers to be twin primes, a result
reminiscent of Clement's theorem but involving indivisibility.

\end{abstract}

\section{Introduction.}

There are some theorems in number theory that correlate the primality of a
number with a congruence. One well-known example is:

\bigskip

\noindent\textbf{Wilson's Theorem:} \emph{A necessary and sufficient condition
for a positive integer }$n\geq2$\emph{ to be a prime is}%
\begin{equation}
(n-1)!+1\equiv0\quad\text{mod }n. \label{Wilson}%
\end{equation}

\noindent In this theorem the primality of a number is connected to the
\textit{divisibility} of the quantity $(n-1)!+1$ by that same number. Here we
first present a result, Theorem 1, similar to the one above, but that differs
from it in an important respect: the primality of the number rests on the
\textit{indivisibility} of a certain quantity by that number. We then present
another result, Theorem 2, that bears on twin primes and follows from Theorem 1.

\section{A theorem correlating primality and indivisibility.}

We have the following theorem:

\bigskip

\noindent\textbf{Theorem 1}. \emph{Let }$m\geq5$\emph{ be a positive integer.
Then a necessary and sufficient condition for it to be a prime is that}%
\begin{equation}
4(m-5)!\equiv\!\!\!\!\!/\;0\quad\text{mod }m.\label{Theorem 1}%
\end{equation}
\emph{There are two exceptions, }$m=6$\emph{ and }$9;$\emph{ the theorem does
not hold for either.}

For example, take $m=8,$ so that $4(8-5)!\equiv24\equiv0$ mod $8,$ so that 8
cannot be prime; and $m=11,$ so that $4(11-5)!\equiv2880\equiv\!\!\!\!\!/\;0$
mod $11,$ so that 11 is a prime. We now prove the theorem.

\noindent\textbf{Proof (necessity). }We first prove necessity. We must show
that if $m$ is prime, it cannot then divide $4(m-5)!.$ Simply notice all
factors in this last expression are smaller than $m;$ since it is prime, it
would not be able to divide $4(m-5)!.$

\noindent\textbf{Proof (sufficiency). }We now prove sufficiency: if
$4(m-5)!\equiv\!\!\!\!\!/\;0$ mod $m,$ then $m$ is a prime.

This can be shown directly for the three smaller values $m=5,7,8$ (for $m=6,9$
the theorem does not hold). Thus, for $m=5$ the congruence $4\cdot
0!\equiv\!\!\!\!\!\!/\;0$ mod $5$ holds, and, in accordance with the
sufficiency of the theorem, $5$ is prime. For $m=7$ the congruence
$4\cdot2!\equiv\!\!\!\!\!/\;0$ mod $7$ holds and 7 is a prime. For $m=8$ we
have instead that $4\cdot3!\equiv0$ mod $8,$ so that 8 should not be a prime.

For $m\geq10$ we proceed by showing the contrapositive to the sufficiency,
that is, if $m$ can be factorized then $4(m-5)!\equiv0$ mod $m$. If $m$ is not
prime then it can be written in the form $m=pq,$ where we take $p$ to be the
smaller of the two factors and $q$ the larger. Observe that, since $p\geq2$,
we have an immediate condition on the maximum value possible for $q$:%
\[
q=m/p\leq m/2.
\]
If this other factor $q$ is present in the expression $4(m-5)!$ then $m=pq$
will divide this expression, since it would certainly contain both its factors
(as long as they do not happen to be equal). This will be the case if
inequality%
\[
m/p\leq m/2\leq m-5
\]
holds, for in this case the factorial $(m-5)!$ would certainly have to contain
the factor $q.$ Since the inequality is satisfied for $m\geq10,$ then both
factors of $m$ must be present in $4(m-5)!,$ and we have proved that
$4(m-5)!\equiv0$ mod $m$.

The above argument does not work in the special case in which $m$ is a squared
prime, but then it is easy to see that $m=$ $p^{2}$ divides $4(m-5)!$ in this
case, too. Simply notice that for $p\geq5$ it happens that $(p^{2})!$ always
contains $p+1$ factors of $p,$ so that $4(p^{2}-5)!$ contains $p-1$ factors of
$p$, which is more than the two needed for $4(p^{2}-5)!$ to be divisible by
$m=p^{2}.\blacksquare$

\section{A new test for twin primes.}

Twin primes are two prime number that differ by two, like 5 and 7. There is a
theorem [\textbf{1}] that gives a test for two numbers being twin primes:

\bigskip

\noindent\textbf{Clement's Theorem: }\emph{A necessary and sufficient
condition for the integers }$n$\emph{ and }$n+2,$\emph{ where }$n\geq2,$\emph{
to be twin primes, is that the congruence }%
\begin{equation}
4\left(  (n-1)!+1\right)  +n\equiv0\quad\text{mod }n(n+2) \label{Clement}%
\end{equation}
\emph{holds true.\noindent}

We can use Theorem 1 to establish another test for twin primes, similar to
Clement's Theorem, except that it involves indivisibility instead of divisibility:

\bigskip

\noindent\textbf{Theorem 2.}\emph{ Take a positive integer }$n\geq3.$\emph{
Then a necessary and sufficient condition for the integers }$n$\emph{ and
}$n+2$\emph{ to be twin primes is that the quantity }%
\begin{equation}
4(n-3)!+2+n \label{Theorem 2}%
\end{equation}
\emph{is divisible by }$n$\emph{ but not by }$n+2$\emph{. There is one
exception to the theorem: }$n=7.$

To illustrate Theorem 2 with only one prime (and not a twin) take $n=13,$ so
that (\ref{Theorem 2}) takes the value 14515215 which divides both 13 and 15.
Since 13 \emph{is} a prime, it must that its supposed twin, 15, is not a
prime, as indeed it is not.

To illustrate Theorem 2 for actual twin primes take $n=5.$ Then expression
(\ref{Theorem 2}) takes the value 15, which is divisible by 5 but not by $7$.
We conclude that 5 and 7 must be twin primes, as indeed they are.

\noindent\textbf{Proof (necessity). }We must show that if $n$ and $n+2$ are
prime, then the congruences $4(n-3)!+2+n\equiv0$ $\operatorname{mod}$ $n$ and
$4(n-3)!+2+n\equiv\!\!\!\!\!/\;0$ $\operatorname{mod}$ $n+2$ are true.

Since $n$ is prime by hypothesis, we know from Wilson's Theorem that
$(n-1)!+1\equiv0$ mod $n.$ We now take $n\geq3,$ and notice that since
$(n-1)!=(n-1)(n-2)[(n-3)!]$ and $(n-1)(n-2)\equiv2$ mod $n,$ Wilson's relation
leads to $2(n-3)!+1\equiv0$ $\operatorname{mod}$ $n.$ In order to be able to
use Theorem 1, we wish to have 4 as a factor of $(n-3)!.$ By multiplying both
sides of the previous congruence by 2 we obtain $4(n-3)!+2\equiv0$
$\operatorname{mod}$ $n,$ or, what is the same,
\begin{equation}
4(n-3)!+n+2\equiv0\quad\operatorname*{mod}\text{ }n. \label{WT2}%
\end{equation}

Let us now rewrite Theorem 1 using the new variable $n=m-2.$ Then the theorem
states that for $n\geq3,$ a necessary and sufficient condition for $n+2$ to be
prime is that $4(n-3)!\equiv\!\!\!\!\!/\;0$ mod $n+2.$ The two exceptions
$m=6$ and 9 become $n=4$ and $7$. Actually $n=4$ is not an exception any more,
since its supposed twin partner, 6, is not a prime either, so that, while
Theorem 1 has two exceptions, Theorem 2 has only one.

Since $n+2$ is prime by hypothesis we can infer from Theorem 1 that%
\begin{equation}
4(n-3)!\equiv4(n-3)!+n+2\equiv\!\!\!\!\!/\;0\quad\operatorname*{mod}\text{
}n+2,\label{T2T1}%
\end{equation}
which is the other congruence we needed for necessity. This concludes the
necessity proof.

\noindent\textbf{Proof (sufficiency). }We must now prove that if the
expression $4(n-3)!+n+2$ is divisible by $n$ but not by $n+2,$ then both $n$
and $n+2$ are prime, with $n\geq3.$

Let us assume that the expression is divisible by $n,$ that is, that%
\[
4(n-3)!+n+2\equiv0\quad\operatorname*{mod}\text{ }n
\]
Then we can retrace some of the steps of (\ref{WT2}) back, except that we
cannot simply divide by 2, since that is sometimes invalid. We are thus left
with the congruence%
\begin{equation}
2(n-1)!\equiv-2\quad\operatorname*{mod}\text{ }n.\label{con}%
\end{equation}
If we assume that $n$ is odd, then we can divide (\ref{con}) by 2 and conclude
that $n$ is prime by Wilson's Theorem. The alternative, to assume that $n$ is
even, leads to contradiction: In this case $n$ can be written in the form
$n=2q,$ and then (\ref{con}) takes the form%
\[
(2q-1)!\equiv-1\quad\operatorname*{mod}\text{ }q.
\]
As the left side of the congruence is divisible by $q,$ we are left with the
contradictory congruence $0\equiv-1$ $\operatorname*{mod}$ $q$. Thus $n\geq3$
has to be odd and prime, or, simply, prime.

Going back to the beginning of the sufficiency proof, let us now assume that%
\[
4(n-3)!+n+2\equiv\!\!\!\!\!/\;0\quad\operatorname*{mod}\text{ }n+2.
\]
We can retrace the steps of (\ref{T2T1}) back and immediately obtain%
\[
4(n-3)!\equiv\!\!\!\!\!/\;0\quad\operatorname*{mod}\text{ }n+2.
\]
Then by Theorem $1$ we conclude that $n+2$ is also prime. This concludes the
sufficiency proof.$\blacksquare$

\bigskip

\noindent REFERENCES

\begin{enumerate}
\item P. A. Clement. Congruences for Sets of Primes, \emph{Amer. Math.
Monthly} \textbf{56} (1949) 23.
\end{enumerate}

\end{document}